\documentclass{article}
\usepackage{epsfig}
\usepackage{amssymb}

\newcommand{\Pol}{{\mbox{\textsf{Pol}}}}
\newcommand{\Tr}{{\mbox{\textsf{Tr}}}}
\newcommand{\kp}{{\Pi}}
\newcommand{\R}{{\mathbb R}}
\newcommand{\C}{{\mathbb C}}
\newcommand{\Z}{{\mathbb Z}}
\newcommand{\HH}{{\mathcal H}}
\renewcommand{\L}{{\mbox{\textsf{Lin}}}}
\newcommand{\sn}{{\ast^{}_\hbar}}
\newcommand{\sg}{{\ast^G_\hbar}}
\newcommand{\sk}{{\ast^K_\hbar}}
\newcommand{\g}{{\mathfrak g}}
\newcommand{\symg}{{\mathcal{S}(\mathfrak{g})}}
\newcommand{\ug}{{\mathcal{U}(\mathfrak{g})}}
\newcommand{\B}{\mathcal{B}}

\newtheorem{theorem}{Theorem}
\newtheorem{definition}{Definition}
\newtheorem{lemma}{Lemma}
\newtheorem{remark}{Remark}
\newtheorem{corollary}{Corollary}
\newenvironment{proof}{\noindent{\it Proof.\/}}{\vskip3mm}
\newenvironment{acknowledgement}{{\vskip3mm}
{\noindent{\bf Acknowledgements.\/}}}{\vskip3mm}

\begin{document}
\textbf{\large Kontsevich star-product on the dual of a Lie algebra}
\vskip10mm
\noindent{Giuseppe Dito}\\
{Laboratoire Gevrey de Math\'ematique physique\\ 
Universit\'e de Bourgogne, BP 47870, F-21078 Dijon Cedex, France}\\
e-mail: \texttt{ditog@u-bourgogne.fr}\vskip10mm
\begin{center}
\textit{Dedicated to the memory of Mosh\'e Flato.}
\end{center}
\vskip10mm
\begin{abstract}
We show that on the dual
of a Lie algebra $\g$ of dimension $d$, the 
star-product recently introduced by M.~Kontsevich is
equivalent to the Gutt star-product on~$\g^\ast$. We give
an explicit expression for the operator realizing 
the equivalence between these star-products.
\end{abstract}

\section{Introduction} 
The study of formal (1-differentiable) deformations of the Lie-Poisson
algebra of functions on a symplectic manifold was initiated in a paper 
by M.~Flato, A.~Lichnerowicz and D.~Sternheimer~\cite{FLS75a}. 
Shortly afterwards, this program was extended
to star-products, i.e., associative deformation of the usual product of
functions, on a symplectic manifold \cite{BFFLS78a} giving,
among others, a profound interpretation of Quantum Mechanics as a deformation
of Classical Mechanics in the direction of the Poisson bracket. 

The existence problem of star-products has been solved by successive steps
from special classes of symplectic manifolds to general Poisson manifolds.
The existence of star-products on any finite dimensional symplectic
manifold was first shown by M.~De~Wilde and P.~Lecomte~\cite{DWLe83b}.
Since then, more geometric proofs have appeared \cite{OMY91,Fedb94a}, and
a proof of existence for regular Poisson manifolds was published
by M.~Masmoudi~\cite{Masm92a}.

For the non-regular Poisson case, first examples of star-products appeared
in \cite{BFFLS78a} in relation with the quantization of angular 
momentum. They were defined on the dual of $\mathfrak{so}(n)$ endowed with its
natural Kirillov-Poisson structure. The case for any Lie algebra follows from 
the construction given by S.~Gutt \cite{Guts83b}
of a star-product on the cotangent
bundle of a Lie group $G$. This star-product restricts to a star-product on
the dual of the Lie algebra of $G$. It translates the associative structure
of the enveloping algebra in terms of functions on the dual of the Lie
 algebra of $G$.

The problem of existence of star-products on any finite dimensional
Poisson manifold was given a solution by M.~Kontsevich \cite{Konm97b}. The
proof is based on an explicit expression of a star-product on
$\R^d$ endowed with a general Poisson bracket, which itself follows from
more general formulae which allowed him to show his formality conjecture
\cite{Konm97a} for $\R^d$ and then for any finite dimensional manifold $M$.
Recently, by using different techniques, 
D.~Tamarkin \cite{Tamd99a} has indicated
another proof of the formality conjecture. This is one of the ingredients
in the most recent fundamental paper by M.~Kontsevich~\cite{Konm99}.

On the dual of a Lie algebra, we have a priori two different star-products: 
Gutt and Kontsevich star-products. D.~Arnal \cite{Arnd98} showed that when the
Lie algebra is nilpotent, these two star-products do coincide. Here we shall
give an elementary proof that in the general case Gutt and Kontsevich 
star-products are equivalent and explicitly construct the equivalence between 
them. For that purpose, we use the notion of Weyl star-products on $\R^d$.
These are star-products enjoying  the following property: $X\sn\cdots\sn X$ 
($k$ factors) is equal to $X^k$  (usual product) for any  linear polynomial 
$X$on $\R^d$ and $k\geq0$. Any star-product on $\R^d$ is equivalent to a 
Weyl star-product and, in the case of the dual of a Lie algebra, Gutt 
star-product is the unique covariant Weyl star-product. From this fact, one
immediately obtains that Kontsevich and Gutt star-products are equivalent.
The equivalence operator is obtained by applying a method used in \cite{Ditg98}
in the context of generalized Abelian deformations.
It turns out that the equivalence operator is an exponential of constant
coefficients linear operators given by the trace of powers of the adjoint map
of the Lie algebra. The formula we have found is closely related to the 
discussion given in \cite{Konm97b} about Lie algebras.
 
The paper is organized as follows. We review Kontsevich construction 
in Section~2. Then we proceed
to the study of Weyl star-products on $\R^d$ and get a characterization
of Gutt star-product. The main results about equivalence is proved 
in Section~4. We end the paper with some remarks on the general Poisson case.

Since the first version of this paper was completed, several preprints dealing
with Kontsevich star-product on the dual of a Lie algebra have appeared 
\cite{ABM99,Katv98,Shob99}. The equivalence result found here
has also been obtained by D.~Arnal, N.~Ben Amar, and M.~Masmoudi \cite{ABM99}
in a completely different approach involving cohomology of Kontsevich graphs.

\section{Kontsevich star-product}
The reader is referred to \cite{BFFLS78a} for 
the general theory on star-products and to the excellent review
by D.~Sternheimer \cite{Sted98} for further details 
and recent applications.

We shall briefly review the construction of a star-product on $\R^d$ given in 
\cite{Konm97b}. Consider $\R^d$ endowed with a Poisson bracket $\pi$. 
We denote by  $(x^1,\ldots,x^d)$ the coordinate system on $\R^d$, the Poisson 
bracket of two smooth functions  $f,g$ is given by 
$\pi(f,g)= \sum_{1\leq i,j\leq n}\pi^{ij}\partial_i f \partial_j g$, 
where $\partial_k$ denotes the 
partial derivative with respect to $x^k$. What follows remain valid if, 
instead of the whole of  $\R^d$, one considers an open subset of it. 
We slightly depart from the notations used in \cite{Konm97b}.

The formula for Kontsevich star-product is conveniently defined
by considering, for each $n\geq 0$, a family of oriented graphs $G_n$. To a 
graph $\Gamma \in G_n$ is associated a bidifferential operator $\B_\Gamma$ and
a weight $w(\Gamma)\in \R$. 
The sum $\sum_{\Gamma\in G_n}
w(\Gamma) \B_\Gamma$ gives us the term at order~$\hbar^n$, i.e., the cochain 
$C_n$ of the star-product. Here is the formal definition of $G_n$.

An oriented graph $\Gamma$ belongs to $G_n$, $n\geq 0$, if:
\begin{itemize}
\item[i)] $\Gamma$ has $n+2$ vertices labeled $\{1,2,\ldots,n,L,R\}$ where $L$
and $R$ stand for Left and Right, respectively, and 
$\Gamma$ has $2n$ oriented edges labeled $\{i_1,j_1,i_2,j_2,\ldots,i_n,j_n\}$;
\item[ii)]  The pair of edges $\{i_k,j_k\}$, $1\leq k \leq n$, starts at 
vertex~$k$;
\item[iii)] $\Gamma$ has no loop (edge starting at some vertex and ending at
that vertex) and no parallel multiple edges (edges sharing the same starting
and ending vertices).
\end{itemize}
When it is needed to make explicit at which  vertex 
$v\in\{1,\ldots,n,L,R\}$ some edge, e.g. $j_k$, is ending at,
we shall use the notation  $j_k(v)$.  

The set of graphs in $G_n$ is finite. For $n\geq1$, the first edge $i_k$ 
starting at vertex $k$ has $n+1$ possible ending vertices 
(since there is no loop), while the second edge $j_k$ has only $n$ possible 
ending vertices, since there is no parallel multiple edges. Thus there are
$n(n+1)$ ways to draw the pair of edges starting at some vertex and therefore 
$G_n$ has $(n(n+1))^n$ elements. For $n=0$, $G_0$ has only one element: 
The graph having as set of vertices $\{L,R\}$ and no edges.

A bidifferential operator $(f,g)\mapsto \B_\Gamma(f,g)$,
$f,g\in C^\infty(\R^d)$,
is associated to each graph $\Gamma\in G_n$, $n\geq1$. 
To each vertex~$k$, $1\leq k\leq n$, one associates the components
$\pi^{i_kj_k}$ of the Poisson tensor, $f$ is associated to the vertex~$L$ and 
$g$ to the vertex $R$. Each edge, e.g. $i_k(v)$ acts by partial differentiation
with respect to $x^{i_k}$ on its ending vertex~$v$. 
There is no better way than to draw the graph $\Gamma$
to illustrate the correspondence 
$\Gamma\mapsto \B_\Gamma$. See \cite{Konm97b}
for a general formula. The graph in Fig.~1 gives the bidifferential operator
$$
\B_\Gamma(f,g) = \sum_{0\leq i_\ast,j_\ast\leq d}
\pi^{i_1j_1}\,\partial_{j_1j_3}\pi^{i_2j_2}\,
\partial_{i_2}\pi^{i_3j_3}\,\partial_{i_1j_2}f\,\partial_{i_3}g.
$$
Notice that for $n=0$, we simply have the usual product of $f$ and $g$.
\begin{figure}
\centerline{\epsfig{file=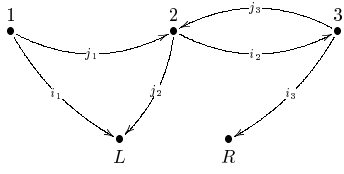,width=45mm}}
\caption{A typical graph in $G_3$.}
\end{figure}

Now let us describe how the weight $w(\Gamma)$ of a graph $\Gamma$ is defined. 
Again the reader is referred to \cite{Konm97b} for details and a nice 
geometrical interpretation of what follows. Let 
$\HH=\{z\in\C\,|\,\mathrm{Im}(z)>0\}$ be the upper half-plane. 
$\HH_n$ will denote the configuration space 
$\{z_1,\ldots,z_n\in\HH\,|\,z_i \neq z_j\ \mathrm{for}\ i\neq j\}$.
$\HH_n$ is an open submanifold of $\C^n$. Let 
$\phi\colon \HH_2\rightarrow \R/2\pi \Z$ be the function:
\begin{equation}\label{phih}
\phi(z_1,z_2)=\frac{1}{2\sqrt{-1}}\mathrm{Log}
\Bigl(\frac{(z_2-z_1)(\bar{z}_2-z_1)}
{(z_2-\bar{z}_1)(\bar{z}_2-\bar{z}_1)}\Bigr).
\end{equation}
$\phi(z_1,z_2)$ is extended by continuity for $z_1,z_2\in \R$, $z_1\neq z_2$.

For a graph $\Gamma\in G_n$, the vertex $k$, $1\leq k\leq n$, is associated
with the variable  $z_k\in\HH$, the vertex $L$ with $0\in\R$, and the vertex 
$R$ with $1\in\R$. 

The weight  $w(\Gamma)$  is defined by integrating an $2n$-form over $\HH_n$:
\begin{equation}\label{w}
w(\Gamma)=\frac{1}{n!(2\pi)^{2n}}\int_{\HH_n}\bigwedge_{1\leq k\leq n}
\Bigl(d\phi(z_k,I_k)\wedge d\phi(z_k,J_k)\Bigr),
\end{equation}
where $I_k$ (resp. $J_k$) denotes the variable or real number associated with
the ending vertex of the edge $i_k$ (resp. $j_k$). For example, the weight
of the graph in Fig.~1 consists in integrating the $6$-form
$d\phi(z_1,0)\wedge d\phi(z_1,z_2)
\wedge d\phi(z_2,z_3)\wedge d\phi(z_2,0)
\wedge d\phi(z_3,1)\wedge d\phi(z_3,z_2)$ on $\HH_3$.
It is clear from the definition of the weights that they are universal
in the sense that they do not depend on the Poisson structure or the
dimension~$d$. 

The origin of the weights has been elucidated by A.~S.~Cattaneo and 
G.~Felder \cite{CaFe99}. These authors  have been able
to construct a bosonic topological field theory on a disc whose perturbation series
(after a finite renormalization taking care of tadpoles) makes 
Kontsevich graphs and weights appear explicitly.

\begin{figure}
\centerline{\epsfig{file=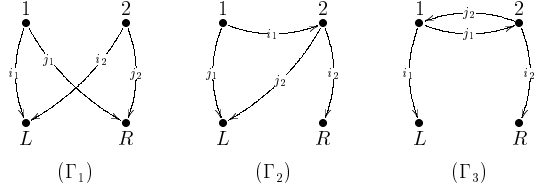,width=70mm}}
\caption{Graphs contributing to $C_2^K$.}
\end{figure}
It is showed in \cite{Konm97b} that the integral in Eq.~(\ref{w})
is absolutely convergent. A pillar result in \cite{Konm97b} is
\begin{theorem}[Kontsevich]\label{konsi}
For any Poisson structure $\pi$ on $\R^d$, the map 
$$(f,g)\mapsto \sum_{n\geq0}\hbar^n
\sum_{\Gamma\in G_n} w(\Gamma)\B_{\Gamma}(f,g),\quad f,g\in C^\infty(\R^d),
$$
defines an associative product.
\end{theorem}

We call this product the Kontsevich star-product and it will be denoted by
$\sk$ and the corresponding cochains  by $C^K_r$. 
Actually the preceding theorem holds if $\pi$ is replaced by any formal 
Poisson bracket $\pi_\hbar= \pi +\sum_{r\geq1}\hbar^r \pi_r$.
Moreover, equivalence classes of star-products are in one-to-one correspondence
with equivalence classes of formal Poisson brackets \cite{Konm97b}.

At first it may seem that any computation involving the graphs becomes rapidly
cumbersome as $\#G_n=(n(n+1))^n$. But the situation is not that bad, 
there are many isomorphic graphs obtained
by permuting the vertices or interchanging the edges 
$\{i_k,j_k\}\rightarrow \{j_k,i_k\}$. These operations do not affect 
$w(\Gamma)\B_\Gamma$ as each factor picks up a minus sign. 
Also in $G_n$, $n\geq2$, there are ``bad''
graphs that can be eliminated right away. These graphs are those for which
the vertices $L$ or $R$ (or both) do not receive any edge. As it should, the
weights associated to these graphs vanish. 

We will illustrate that by giving the explicit form of the second cochain 
$C_2$ which requires at the end the computation of only $3$ graphs 
(notice that $\# G_2=36$).

The graphs in Fig.~2 have weights $w(\Gamma_1)=1/8$, 
$w(\Gamma_2)=1/24$, $w(\Gamma_3)= -1/48 $. By counting the symmetries, 
the graph $\Gamma_1$ contributes $4$ times, $\Gamma_2$ contributes $8$ times, 
and  $\Gamma_3$ contributes $8$ times. There is also a sister-graph for 
$\Gamma_2$ which is obtained by performing $L\leftrightarrow R$ which 
contributes also $8$ times. Taking into account that there are $8$ 
``bad'' graphs, we a have a  total of $36$ graphs, and we find that:
\begin{eqnarray}\label{c2}
C^K_2(f,g)& = &\frac{1}{2}\,\pi^{i_1j_1}\,\pi^{i_2j_2}\,\partial_{i_1 i_2} f\, 
\partial_{j_1 j_2} g\nonumber\\
&&+ \frac{1}{3}\,\pi^{i_1j_1}\,\partial_{i_1}\pi^{i_2j_2}\,
(\partial_{j_1 j_2} f\,\partial_{i_2} g +\partial_{i_2} f\, 
\partial_{j_1 j_2}g)\\
&&-\frac{1}{6}\,\partial_{j_2}\pi^{i_1j_1}\,
\partial_{j_1}\pi^{i_2j_2}\,\partial_{i_1} f\, \partial_{i_2}g,\nonumber
\end{eqnarray}
where summation over repeated indices is understood.

\section{Weyl star-products on  $\R^d$}
Let $\pi$ be a general Poisson structure on $\R^d$.  Let $\Pol$ be
the algebra of polynomials in the variables $x^{{1}},\ldots, x^{{d}}$ and
let $\L$ denote the subspace of linear homogeneous polynomials.
Let $\sn$ be a star-product on $(\R^d,\pi)$. We shall show that $\sn$ is 
(differentially) equivalent to a 
star-product $\sn'$ having the following property:
\begin{equation}\label{sym}
X^{\sn'k}= X^k,\quad \forall k\geq 0,\forall X\in \L,
\end{equation}
where $X^{\sn'k}=X\sn'\cdots\sn'X$ ($k$ factors). This is reminiscent of
Weyl ordering in Quantum Mechanics and we introduce:
\begin{definition}
A star-product on $(\R^d,\pi)$ satisfying Eq.~(\ref{sym}) 
is called a Weyl star-product.
\end{definition}
The consideration of this kind of star-products amount to 
generalized Abelian deformations \cite{DiFl97a,Ditg98}. We recall the
proof of the following:
\begin{theorem}\label{weyl}
Any star-product on  $(\R^d,\pi)$ is equivalent to a Weyl star-product.
\end{theorem}
\begin{proof}
The  proof consists to establish the differentiability of the following
$\R[[\hbar]]$-linear map 
$\rho^{}\colon\Pol[[\hbar]]\rightarrow C^\infty(\R^d)[[\hbar]]$
uniquely defined by:
\begin{equation}\label{rho}
\rho^{}(X^k)=X^{\sn k}, \quad \forall k\geq 0,\forall X\in \L. 
\end{equation}
The map $\rho^{}$ is a formal sum of linear maps 
$\rho^{}=\sum_{i\geq 0}\hbar^i \rho^{}_i$ with $\rho^{}_0$ being the identity 
map on~$\Pol$. We will show that the $\rho^{}_r$'s are differential operators.
By definition $\rho^{}_r(1)=\rho^{}_r(X)=0$ for $r\geq1$ and $X\in\L$. It is 
easy to see from Eq.~(\ref{rho}) that the $\rho^{}_r$'s satisfy the following 
recurrence relation for  $k\geq1,r\geq1$:
\begin{equation}\label{recrho}
-\delta\rho^{}_r(X,X^{k-1})= C_r(X,X^{k-1})+
\sum_{a+b=r \atop a,b\geq1}C_a(X,\rho^{}_b(X^{k-1})),
\end{equation}
(the $C_r$'s are the 2-cochains of the star-product). For $r=1$ the sum
on the right-hand side is omitted and $\delta$ is the 
Hochschild differential. Before going further we need a lemma.
\begin{lemma}\label{diff}
Let $\psi\colon\Pol\rightarrow C^\infty(\R^d)$ be an $\R$-linear map such that
$\psi(1)=\psi(X)=0$, for $X\in \L$, and let 
$\phi\colon C^\infty(\R^d)\times C^\infty(\R^d)\rightarrow C^\infty(\R^d)$
be a bidifferential operator vanishing on constants. If $\psi$ satisfies
\begin{equation}\label{phi}
\delta\psi(X,X^{k-1})=\phi(X,X^{k-1}),\quad \forall k\geq1, \forall X\in\L,
\end{equation}
then there exists a differential operator $\eta$ on $\R^d$ such that 
$\psi=\eta|_{\Pol}$.
\end{lemma}
\begin{proof} For two functions $f,g$, let  
$\sum_{I,J}\phi^{I,J}\partial_If \partial_J g$ be the expression of $\phi(f,g)$ 
in local coordinates, 
where $I$ and $J$ are multi-indices and $\phi^{I,J}$ is a smooth function
vanishing for $|I|$ or $|J|$ greater than some integer 
(for $I=(i_1,\ldots,i_n)$, $|I|$ denotes its length $i_1+\cdots+i_n$).
In Eq.~(\ref{phi}), only first derivatives  can be applied to the first argument
of $\phi$ and one can check the following series of equalities:
\begin{equation}\label{eta}
\phi(X,X^{k-1})=\sum_{i,J}\phi^{i,J}\partial_iX\partial_J X^{k-1}
={1\over k}\sum_{i,J}\phi^{i,J}\partial_{iJ}X^k = \delta \eta(X,X^{k-1}),
\end{equation}
where $\eta=-\sum_{i,J}{1\over 1+|J|}\phi^{i,J}\partial_{iJ}$. Now 
Eq.~(\ref{phi}) can be written as $\delta(\psi-\eta)(X,X^{k-1})=0$ and by
observing that $\psi-\eta$ vanishes on $1$ and $X$ we get that $\psi-\eta=0$
on $\Pol$. This shows the lemma.
\end{proof}

The term of order $1$ in Eq.~(\ref{recrho}) yields 
$(\delta\rho^{}_1 + C_1)(X,X^{k-1})=0$. We have that $C_1=\pi+\delta\theta_1$ 
for some differentiable $1$-cochain $\theta_1$, which can be chosen such that
$\theta_1(X)=0$ for $X\in \L$ by adding an appropriate Hochschild $1$-cocycle
(i.e., a vector field). Then as before  
$\delta(\rho^{}_1 + \theta_1)(X,X^{k-1})=0$ gives us $\rho^{}_1=-\theta_1$ on
$\Pol$ showing that $\rho^{}_1$ is a differential operator.

With the help of Lemma~\ref{diff}, a simple recurrence on $r$ in 
Eq.~(\ref{recrho}) shows that for each $r\geq1$, $\rho^{}_r$  coincides with
the restriction of a differential operator to $\Pol$. Clearly the map $\rho$
can be naturally extended to an $\R[[\hbar]]$-linear map on  
$C^\infty(\R^d)[[\hbar]]$. We still denote this extension by $\rho$.

The map $\rho$ is invertible as $\rho^{}_0$ is the identity map and we can
use it to define an equivalent star-product $\sn'$ to  $\sn$ by:
\begin{equation}
\rho(f\sn'g)=\rho(f)\sn\rho(g),\quad f,g\in C^\infty(\R^d).
\end{equation}
Notice that $X^{\sn'k}=\rho^{-1}(\rho(X)^{\sn k})=\rho^{-1}(X^{\sn k})=X^k$ for
$\forall k\geq 0$ and $\forall X\in \L$,
therefore $\sn'$ is a Weyl star-product.
\end{proof}

\subsection{Gutt star-product on $\g^\ast$} 
Let $G$ be a real finite-dimensional group of dimension $d$. 
The Lie algebra of $G$ is denoted by $\g$ and its dual by $\g^\ast$. 
The universal enveloping algebra (resp.~symmetric algebra) of $\g$ is denoted
by $\ug$ (resp.~$\symg$). Also we denote by $\Pol(\g^\ast)$ the space of
polynomials on $\g^\ast$.

It is well known that the space of smooth functions on $\g^\ast$ carries a 
natural Poisson structure defined by the Kirillov-Poisson bracket denoted 
by $\kp$. Fix a basis for $\g$, let $C_{ij}^k$ be the structure constants 
in that basis, and let $(x^1,\ldots,x^d)$ be the coordinates
on $\g^\ast$. Then  the Kirillov-Poisson bracket is defined by:
\begin{equation}\label{kp}
\kp(f,g)=\sum_{1\leq i,j, k\leq d} x^k\, C_{ij}^k\, \partial_i\, f \partial_j g,
\quad f,g\in C^\infty(\g^\ast).
\end{equation}
Of course this definition is independent of the chosen basis for $\g$.

S.~Gutt has defined a star-product on the cotangent bundle $T^\ast G$ of a Lie 
group $G$ \cite{Guts83b}. When one restricts this star-product between 
functions not depending on the base point in $G$, one gets a star-product 
on $\g^*$. We shall call the induced product on $C^\infty(\g^\ast)$, 
Gutt star-product  on $\g^\ast$ and denote it by $\sg$. 
Gutt star-product on $\g^\ast$ can also
be directly obtained by transporting the algebraic structure of the enveloping 
algebra $\ug$ of $\g$ to the space of polynomials on 
$\g^\ast$. This is achieved through the natural
isomorphism between $\Pol(\g^\ast)$ and $\symg$ and 
with the help of the symmetrization map $\sigma\colon\symg\rightarrow\ug$.
The product between two homogeneous elements $P$ and $Q$ in 
$\symg\sim\Pol(\g^\ast)$ of degrees $p$ and $q$, respectively, is given by:
\begin{equation}\label{gutt}
P\sg Q = \sum_{0\leq r\leq p+q-1}(2\hbar)^r 
\sigma^{-1}((\sigma(P)\cdot\sigma(Q))_{p+q-r}),
\end{equation}
where $\cdot$ is the product in $\ug$ and, for $v\in \ug$, $(v)_k$ means 
the $k$-th component of $v$ with respect to the associated grading of $\ug$.
Formula~(\ref{gutt}) defines an associative deformation of the usual 
product on $\Pol(\g^\ast)$ which admits a unique extension to  
$C^\infty(\g^\ast)$.

As a  direct consequence of Eq.~(\ref{gutt}) we have that $\sg$ is a Weyl
star-product on $(\g^\ast,\kp)$. Moreover the following property is easily
verified:
$$
X\sn Y - Y\sn X = 2\hbar \kp(X,Y),\quad X,Y\in \L(\g^\ast),
$$
where $\L(\g^\ast)$ is the subspace of homogeneous polynomials of degree $1$
on $\g^\ast$. Star-products on  $(\g^\ast,\kp)$ satisfying the preceding 
relation are called $\g$-covariant star-products. Actually, there is a
characterization of  Gutt star-product given by:
\begin{lemma}\label{equigutt}
Gutt star-product is the unique $\g$-covariant Weyl star-product
on  $(\g^\ast,\kp)$. Any  $\g$-covariant star-product on $(\g^\ast,\kp)$ is
equivalent to Gutt star-product.
\end{lemma}
\begin{proof}
Any star-product $\sn$  on  $(\g^\ast,\kp)$ is determined by the quantities
$\exp(X)\ast \exp(Y)$, $X,Y\in \L(\g^\ast)$. Suppose that $\sn$ is a 
$\g$-covariant Weyl star-product, then the star-exponential of 
$X\in\L(\g^\ast)$ defined by:
$$
\exp_\sn(X)=\sum_{k\geq0}\frac{1}{k!}X^{\sn k},
$$
coincides with the usual exponential $\exp(X)$. The covariance property of 
$\sn$ allows to use the Campbell-Hausdorff formula:
\begin{equation}\label{exp}
\exp(X)\sn \exp(Y)=\exp_\sn(CH_\hbar(X,Y))
=\exp(CH_\hbar (X,Y)),
\end{equation}
where $CH_\hbar (X,Y)$ is the usual Campbell-Hausdorff series with respect
to the bracket $[X,Y]=2\hbar\kp(X,Y)$. For $X,Y\in\L(\g^\ast)$,
$CH_\hbar(X,Y)$ is an element of $\L(\g^\ast)[[\hbar]]$. We still have
$CH_\hbar(X,Y)^{\sn k}=CH_\hbar(X,Y)^k$ for $k\geq0$. It follows from
Eq.~(\ref{exp}) that there is at most one $\g$-covariant Weyl star-product on 
$\g^\ast$, i.e., Gutt star-product.
The second statement of the lemma follows from Theorem~\ref{weyl} and from
the fact that the equivalence operator $\rho$ preserves the covariance property,
i.e., $\rho(X)=X$ for $X\in\L(\g^\ast)$ (cf. Eq.~(\ref{rho})).
\end{proof}

From Eq.~(\ref{exp}), one can derive an explicit expression for the cochains
of Gutt star-product. Denote by $c_i$, $i\geq1$, the Campbell-Hausdorff
coefficients: $c_1(X,Y)= X+Y$, $c_2(X,Y)= \frac{1}{2}[X,Y]$, etc. The term
of order $\hbar^r$ in $\exp(X)\sg\exp(Y)$ for $X,Y\in\L(\g^\ast)$ is obtained
by expanding in powers of $\hbar$ the right hand-side in Eq.~(\ref{exp}),
it is given by:
\begin{eqnarray}\label{productexp}
&&C_r^G(\exp(X),\exp(Y))\\
&&\quad=2^r\exp(X+Y)\sum_{1\leq k\leq r}
\sum_{m_1>\cdots>m_k\geq1\atop{ n_1,\ldots n_k \geq 1\atop
m_1 n_1+\cdots+ m_k n_k=r}}
\prod_{1\leq j \leq k}
\frac{1}{n_j!}(c_{m_j+1}(X,Y))^{n_j},\nonumber
\end{eqnarray}
where the bracket $[X,Y]$ in the $c_i$'s is taken to be $\kp(X,Y)$. For $r=2$,
we easily get the differential expression for $C^G_2(f,g)$, 
$f,g\in C^\infty(\g^\ast)$:
$$
\frac{1}{2!}\,\kp^{i_1j_1}\,\kp^{i_2j_2}\,\partial_{i_1 i_2} f\, 
\partial_{j_1 j_2} g
+ \frac{1}{3}\,\kp^{i_1j_1}\,\partial_{i_1}\kp^{i_2j_2}\,
(\partial_{j_1 j_2} f\,\partial_{i_2} g
+\partial_{i_2} f\, \partial_{j_1 j_2}g).
$$
Comparing with the general expression of the second cochain of Kontsevich 
star-product given by (\ref{c2}), we see that in general Kontsevich and 
Gutt star-products differ. Notice that the extra term in (\ref{c2}) is a
Hochschild $2$-coboundary.

It is instructive to derive from  Eq.~(\ref{productexp}) an expression
for $X\sg g$, $X\in \L(\g^\ast)$, $g\in  C^\infty(\g^\ast)$.
Using the standard recurrence formula for the $c_i$'s, it is easy to establish
that
\begin{eqnarray}\label{ci}
&&c_i(0,X) = c_i(X,0) =0,\quad i\geq 2;\nonumber \\
&&\nonumber \\
&&{\partial\over\partial s}c_i(sX,Y)|_{s=0}
= {B_{i-1}\over (i-1)!}(ad_Y)^{i-1}(X),\quad i\geq 2;
\end{eqnarray}
for $X\in\L(\g^\ast)$, $s\in\R$, $ad_Y$ is the adjoint map $X\mapsto[Y,X]$,
and the $B_n$'s are the Bernoulli numbers. The substitution $X\rightarrow
s X$ in Eq.~(\ref{productexp}) and differentiation with respect to $s$
gives:
$$
C_r^G(X,\exp(Y))= \frac{2^r B_r}{r!} (ad_Y)^r(X)\, \exp(Y),
$$
which leads to 
\begin{equation}\label{productxg}
C_r^G(X,g)= (-1)^r\frac{2^r B_r}{r!}\sum_{1\leq i_\ast,j_\ast\leq d}
\kp^{i_1j_1} \partial_{i_1}\kp^{i_2j_2} 
\cdots \partial_{i_{r-1}}\kp^{i_rj_r} \partial_{i_r}X
\partial_{j_1\cdots j_r}g.
\end{equation}
\begin{remark}
From Eq.~(\ref{productexp}) or, better, Eq.~(\ref{productxg}), it is clear
that the weights of the graphs appearing in Gutt star-product are
essentially products of Bernoulli numbers.
\end{remark}

\section{Equivalence}
In this section, as in Sect.~3, we consider a Lie algebra $\g$ 
of dimension $d$ and use the notations previously introduced. We have seen
that in general Kontsevich and Gutt star-products are not identical.
We will show that they are equivalent and explicitly determine the equivalence
operator by computing a subfamily of graphs.

\begin{lemma}\label{konsicov}
Kontsevich star-product $\sk$ on $(\g^\ast,\kp)$ 
is a $\g$-covariant star-product.
\end{lemma}
\begin{proof}
We just need to see what kind of graphs contribute to 
$X\sk Y$, $X,Y\in\L(\g^\ast)$.  The graphs for $C^K_r(X,Y)$ must
be such that the vertices  $L$ and $R$ receive only one edge, respectively.
For $r=1$, we simply have the Poisson bracket $\kp$. If $r\geq2$, 
we need to draw $2r-2$ edges in such a way that each vertex $k$, $1\leq k\leq r$
receives at most one edge (since the Poisson bracket $\kp$ is linear 
in the coordinates) and this is possible only if $2r-2\leq r$, 
i.e., $r\leq2$. For $r=2$, the only
graph contributing (up to symmetry factors) is the graph $\Gamma_3$ in 
Fig.~2, whose  associated bidifferential operator $\B_{{\Gamma_3}}$
is symmetric. Thus we have  $X\sn Y - Y\sn X = 2\hbar \kp(X,Y)$. 
\end{proof}
As a consequence of Lemmas~\ref{equigutt} and \ref{konsicov} we have:
\begin{corollary}
On the dual of a Lie algebra, Kontsevich and Gutt star-products are equivalent.
\end{corollary}

The formal series of differential operators realizing the equivalence 
between Kontsevich and Gutt star-products 
is the map $\rho$ defined in the proof of Theorem~\ref{weyl}. We have
$\rho(f\sg g)=\rho(f)\sk \rho(g)$, $\forall f,g\in C^\infty(\g^\ast)$ and
in the present situation $\rho$ is defined by $\rho(X^k)= X^{\sk k}$, 
$X\in\L(\g^\ast)$, $k\geq0$.
\begin{figure}
\centerline{\epsfig{file=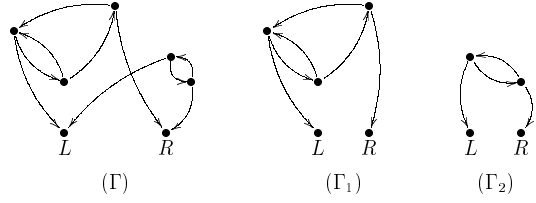,width=70mm}}
\caption{Union of graphs: $\Gamma=\Gamma_1\cup \Gamma_2$.}
\end{figure}
We will see (cf. Theorem~\ref{equiva}) that to solve the
recurrence relation (\ref{recrho}) satisfied by the $\rho_r$'s,
it is sufficient to consider graphs contributing to $C^K_r(X,X^{k})$.
The graphs
having a non-trivial contribution must have only one edge
ending at vertex $L$, e.g. $i_k$, and the other edge $j_k$ must end
at some vertex $k'\neq k$, $1\leq k'\leq r$. 

We shall say that a graph $\Gamma\in G_r$ is the union of two subgraphs 
$\Gamma_1\in G_{r_1}$ and $\Gamma_2\in G_{r_2}$ with  $r_1 + r_2=r$,
if the subset  $(1,\ldots,r)$ of the set of vertices of $\Gamma$ can
be split into two parts $(a_1,\ldots,a_{r_1})$ and $(b_1,\ldots,b_{r_2})$ 
such that there is no edge between these two subsets of vertices. A graph 
that is not the union of two subgraphs is called indecomposable.

By recalling the definition of the weight of a graph, 
the following is straightforward:
\begin{lemma}\label{wunion}
If a graph $\Gamma\in G_r$ is the union of two 
subgraphs $\Gamma_1$ and $\Gamma_2$, respectively, in $G_{r_1}$ and
in $G_{r_2}$ with $r_1 + r_2=r$, then $w(\Gamma)=\frac{r_1!r_2!}{r!}
w(\Gamma_1)w(\Gamma_2)$.
\end{lemma}
\begin{figure}
\centerline{\epsfig{file=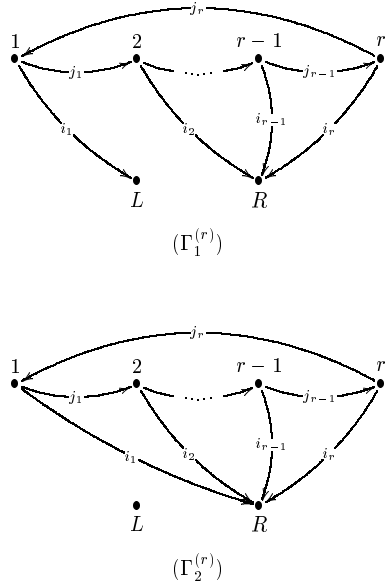,width=50mm}}
\caption{Indecomposable graphs.}
\end{figure}
In view of this lemma, we just need to confine ourselves to indecomposable 
graphs whose union is contributing to $C^K_r(X,X^{k})$.
\begin{lemma}\label{inde}
For $r\geq2$, up to an isomorphism of graphs, 
an indecomposable graph in $G_r$  contributing to $C^K_r(X,X^{k})$ 
falls into one of the two types illustrated in Fig.~4.
\end{lemma}
\begin{proof}
As the vertex $L$ can receive at most one edge, we distinguish two cases.

\noindent\textit{i) The vertex $L$ receives no edge.} We will see that the 
vertex $R$ must receive exactly $r$ edges. If there are strictly more than 
$r$ edges ending at vertex $R$, then there must be a vertex $k$, 
$1\leq k\leq r$, such that the edges $i_k$ and $j_k$ are ending at vertex 
$R$. This is excluded by definition of $G_r$.
If there are strictly less than $r$ edges ending at vertex $R$, 
then at least one of the vertices $(1,\ldots,r)$ must receive two or more edges 
and the bidifferential operator associated to such a graph is vanishing since 
the Poisson bracket is a linear function of the coordinates. We are left with 
the case where exactly $r$ edges are ending at vertex $R$. Then every vertex 
in  $(1,\ldots,r)$ must receive exactly one edge and, up to an isomorphism, 
there is precisely one such a graph, i.e., graph $\Gamma^{(r)}_2$ in Fig.~4.

\noindent\textit{ii) The vertex $L$ receives one edge.} For this case,
the vertex $R$ receives $r-1$ edges.  
By relabeling the vertices and the edges, we may suppose that the edge 
ending at vertex $L$ is $i_1$. Then the second edge
starting at the vertex $1$, i.e., $j_1$, cannot end at vertex $R$ because, by
skew-symmetry of the Poisson bracket, the associated bidifferential operator 
is vanishing on $(X,X^{k})$. Hence we may suppose that the edge $j_1$ is
ending at vertex $2$. We still have $2r-2$ edges starting from the vertices
$(2,\ldots,r)$ to draw. Let $a$ be the number of edges ending at vertex $R$
and let $b$ be the total number of edges ending at the vertices 
$(1,\ldots,r)$. We have $2r-1=a+b$. Since each vertex in $(1,\ldots,r)$
can receive at most one edge, we have 
that $b\leq r$ and it follows that $a\geq r-1$.
If $a> r-1$, it means that there are  parallel multiple edges
between at least one of the vertices $(2,\ldots,r)$ and the vertex $R$. Hence
the vertex $R$ must receive $r-1$ edges. Clearly every such edge must start
at one the vertices $(2,\ldots,r)$. The other $r-1$ edges must
end at the vertices  $(1,3,\ldots,r)$. Thus, up to an isomorphism, we find
that there is only the  graph $\Gamma^{(r)}_1$ in Fig.~4 for this case.
\end{proof}
The preceding lemma tells us that graphs contributing to
$C^K_r(X,X^{k})$ must be of the form:
$$
 \Gamma^{(a)}_1\cup\Gamma^{(b_1)}_2\cup\cdots\cup\Gamma^{(b_s)}_2,
$$
with $a+b_1+\cdots+b_s=r$. Notice that there can be only one graph of
the type $\Gamma^{(a)}_1$, since the vertex $L$ can receive only one edge.
Quite a bit of simplification is allowed by
\begin{lemma}\label{w2}
For $r\geq2$, the weight of the graph $\Gamma^{(r)}_2$ in Fig.~4 vanishes.
\end{lemma}
\begin{proof}
The form $\bigwedge_{1\leq k\leq r}d\phi(z_k,1)\wedge d\phi(z_k,z_{k+1})$,
where $z_{r+1}\equiv z_1$, is $0$. This easily follows from a simple recurrence
using explicit expressions for the forms $d\phi(z_i,z_j)$.
\end{proof}
When they appear alone, the graphs $\Gamma^{(r)}_2$ constitute an example 
of what was called ``bad'' graphs in Sect.~2.

\begin{lemma}\label{op2}
For $r\geq2$,  up to an isomorphism,
the only graph contributing to $C^K_r(X,X^{k})$ is the graph $\Gamma^{(r)}_1$ 
in Fig.~4. The associated bidifferential operator has constant coefficients
and is given by:
\begin{equation}\label{trace}
\B_{\Gamma^{(r)}_1}(f,g)= \sum_{1\leq i_\ast\leq d}
\Tr(ad_{x^{i_1}}\cdots ad_{x^{i_r}})\,
\partial_{i_1}f\,\partial_{i_2\cdots i_r}g,
\quad  r\geq 2.
\end{equation}
\end{lemma}
\begin{proof}
The first statement follows directly from Lemmas~\ref{wunion}, \ref{inde}, 
and~\ref{w2}.
The bidifferential operator for  the graph $\Gamma^{(r)}_1$
is $$
\B_{\Gamma^{(r)}_1}(f,g)=\sum_{1\leq i_\ast,j_\ast\leq d} 
\partial_{j_r}\kp^{i_1j_1}\partial_{j_1}\kp^{i_2j_2}\cdots
\partial_{j_{r-1}}\kp^{i_rj_r}\,\partial_{i_1}f\,\partial_{i_2\cdots i_r}g,
$$
clearly it has constant coefficients and using the expression~(\ref{kp}) 
for $\kp$, we see that the previous equation can
be written as a trace of adjoint maps.
\end{proof}
The computation of the weights of the graphs $\Gamma^{(r)}_1$
is a delicate question. The presence of cycles (wheels) does not
allow do derive a simple recurrence relation among the weights.
A direct calculation for $\Gamma^{(2)}_1$ using residues 
gives a weight equals to $-1/48$, but this method becomes
unpractical for $r\geq3$.

The isomorphic graphs to $\Gamma^{(r)}_1$ are obtained by permuting 
the vertices $(1,\ldots,r)$ and alternating the edges 
$\{i_k,j_k\}\rightarrow \{j_k,i_k\}$ for $1\leq k\leq r$, 
thus we get a symmetry factor $r!2^r$. Hence we have
$C^K_r(X,X^{k})= 2^r r! w(\Gamma^{(r)}_1) 
\B_{\Gamma^{(r)}_1}(X,X^{k})$.

\begin{theorem}\label{equiva}
On any finite-dimensional Lie algebra, the equivalence $\rho$ between
Kontsevich and Gutt star-products: $\rho(f\sg g)=\rho(f)\sk \rho(g)$,  
is given by
\begin{equation}\label{rhoe}
\rho =\exp\Bigl(\sum_{r\geq2} \hbar^r 2^r (r-1)!  w(\Gamma^{(r)}_1) D_r\Bigr),
\end{equation}
where the $D_r$'s are differential operators with constant coefficients:
$$
D_r=\sum_{1\leq i_\ast\leq d} \Tr(ad_{x^{i_1}}\cdots ad_{x^{i_r}})\,
\partial_{i_1\cdots i_r}.
$$
\end{theorem}
\begin{proof}
Recall that $\rho$ is defined by $\rho(X^k)=X^{\sk k}$, $X\in\L(\g^\ast)$,
$k\geq0$. It was shown that the $\rho^{}_r$'s in 
$\rho=I +\sum_{r\geq1} \hbar^r \rho^{}_r$ are differential operators.
Here we have only to solve the recurrence relation for
$\rho_r$ appearing in the proof of Theorem~\ref{weyl}:
\begin{equation}\label{recrhob}
-\delta\rho^{}_r(X,X^{k-1})= C^K_r(X,X^{k-1})+
\sum_{a+b=r \atop a,b\geq1}C^K_a(X,\rho^{}_b(X^{k-1})),
\end{equation}
where $X\in\L(\g^\ast),k\geq1,r\geq1$.
According to Lemma~\ref{diff}, there exist differential
operators $\eta_r$ such that $C^K_r(X,X^{k})=\delta\eta^{}_r (X,X^{k})$.
From Lemma~\ref{op2} it follows that
$$
\eta^{}_r = - 2^r (r-1)!  w(\Gamma^{(r)}_1)\sum_{1\leq i_\ast\leq d}
\Tr(ad_{x^{i_1}}\cdots ad_{x^{i_r}})\,\partial_{i_1\cdots i_r},\quad r\geq2.
$$ 
For each $r\geq2$, $\eta^{}_r$ is a differential operator with constant 
coefficients  and is homogeneous of degree $r$ in the derivatives.
To a differential operator $\eta$  on $\g^\ast$ with
constant coefficients we can associate a polynomial $\hat \eta$ 
on $\g\sim \L(\g^\ast)$. Here we have  
$\hat\eta^{}_r(X) = - 2^r (r-1)!  w(\Gamma^{(r)}_1)
 \Tr((ad_X)^r)$ and one can check that  
\begin{equation}\label{poly}
\delta\eta^{}_r (X,X^{k})= -\frac{r}{k}\eta^{}_r(X^k)
=-\frac{r}{k}\,\frac{k!}{(k-r)!}\,\hat\eta^{}_r(X)\,X^{k-r}.
\end{equation}

The preceding implies that the $\rho_r$'s, $r\geq1$,
have constant coefficients and are homogeneous of degree $r$. 
We have $\rho^{}_1=0$ and  a
recurrence on $r$ in  Eq.~(\ref{recrhob}) shows the property
for all of the  $\rho_r$'s.

Using Eq.~(\ref{poly}) we can express Eq.~(\ref{recrhob}) in terms of the
polynomials $\hat\rho^{}_r$ and  $\hat\eta^{}_r$ and find that:
$$
\hat\rho^{}_r(X)=-\hat\eta^{}_r(X) -\frac{1}{r}\sum_{a+b=r \atop a,b\geq1}
a \hat\eta^{}_a(X)\hat\rho^{}_b(X).
$$
By defining $\hat\eta^{}_0(X)$ to be identically equal to zero and
$\hat\rho^{}_0(X)$ to be $1$, we can
rewrite the previous equation as
\begin{equation}\label{recpolyb}
r\, \hat\rho^{}_r(X)=-\sum_{a+b=r \atop a,b\geq0}
a \hat\eta^{}_a(X)\hat\rho^{}_b(X),
\end{equation}
then by considering the formal series 
$\hat\rho(X)\equiv I +\sum_{r\geq1} \hbar^r \hat\rho^{}_r(X)$
and $\hat\eta(X)\equiv \sum_{r\geq2} \hbar^r \hat\eta^{}_r(X)$ (recall
that $\eta^{}_1=0$), we see that Eq.~(\ref{recpolyb}) simply states
that $\hat\rho'(X)=-\hat\eta'(X)\hat\rho(X)$ where the prime denotes
formal derivative with respect to $\hbar$. 
Thus $\hat\rho(X)=\exp(-\hat\eta(X))$
and Eq.~(\ref{rhoe}) follows.
\end{proof}
For a nilpotent Lie algebra, all of
the operators $D_r$ do vanish. Hence we deduce the result of \cite{Arnd98}:
\begin{corollary}
For a nilpotent Lie algebra, Kontsevich star-product coincides
with Gutt star-product.
\end{corollary}
\subsection{Remarks}
The equivalence between Kontsevich and Gutt star-products shows us
that, in the linear Poisson case, graphs with cycles play no role
with respect to the associativity of the product. Here the contribution
of these graphs is completely absorbed into the equivalence operator.
In other words: Weights of the graphs $\Gamma^{(r)}_1$ can be chosen
arbitrarily and they do not affect the associativity of the star-product.

We suspect that the situation described above is the general one, i.e.,
for $\mathbb{R}^d$ endowed with any Poisson structure $\pi$, 
it would be possible to get a new star-product by removing graphs
with cycles in Kontsevich's construction. We conjecture that 
the Weyl star-product associated with Kontsevich star-product $\sk$
on $(\mathbb{R}^d,\pi)$ contains no cycle, and it is obtained from
$\sk$ by ignoring the graphs with cycles.
\begin{acknowledgement}
Discussions with M.~Flato and D.~Sternheimer were at the 
origin of this paper and I am the most grateful to both of them
for remarks and encouragement.
Most of the work presented here was done while the author
was visiting RIMS with a JSPS grant,
and it a pleasure to thank Prof. I.~Ojima for warmest hospitality. 
\end{acknowledgement}

\end{document}